\documentclass[reqno,centertags, 12pt]{amsart}
\usepackage{amsmath,amsthm,amscd,amssymb}
\usepackage{latexsym}
\usepackage{graphicx}

\newcommand{\bbC}{{\mathbb{C}}}
\newcommand{\bbD}{{\mathbb{D}}}

\newcommand{\bdone}{{\boldsymbol{1}}}

\newcommand{\calC}{{\mathcal C}}
\newcommand{\calH}{{\mathcal H}}

\newcommand{\calM}{{\mathcal M}}

\newcommand{\lb}{\label}
\newcommand{\f}{\frac}

\newcommand{\ol}{\overline}
\newcommand{\ti}{\tilde  }
\newcommand{\wti}{\widetilde  }

\newcommand{\dist}{\text{\rm{dist}}}
\newcommand{\num}{\text{\rm{Num}}}
\newcommand{\prob}{\text{\rm{Prob}}}
\newcommand{\intt}{\text{\rm{int}}}

\newcommand{\spec}{\text{\rm{spec}}}

\newcommand{\bi}{\bibitem}

\newcommand{\beq}{\begin{equation}}
\newcommand{\eeq}{\end{equation}}
\newcommand{\ba}{\begin{align}}
\newcommand{\ea}{\end{align}}
\newcommand{\veps}{\varepsilon}

\DeclareMathOperator{\Real}{Re}
\DeclareMathOperator{\ran}{Ran}

\DeclareMathOperator*{\slim}{s-lim}

\allowdisplaybreaks
\numberwithin{equation}{section}
\newtheorem{theorem}{Theorem}[section]

\newtheorem*{t1}{Theorem 1}
\newtheorem*{t2}{Theorem 2}
\newtheorem*{t3}{Theorem 3}
\newtheorem*{t4}{Theorem 4}
\newtheorem{proposition}[theorem]{Proposition}
\newtheorem{lemma}[theorem]{Lemma}
\newtheorem{corollary}[theorem]{Corollary}
\theoremstyle{definition}

\theoremstyle{remark}
\newtheorem*{remark}{Remark}
\newtheorem*{remarks}{Remarks}
\theoremstyle{definition}
\newtheorem*{definition}{Definition}
\newcommand{\abs}[1]{\lvert#1\rvert}

\newcounter{smalllist}
\newenvironment{SL}{\begin{list}{{\rm\roman{smalllist})}}{%
\setlength{\topsep}{0mm}\setlength{\parsep}{0mm}\setlength{\itemsep}{0mm}%
\setlength{\labelwidth}{2em}\setlength{\leftmargin}{2em}\usecounter{smalllist}%
}}{\end{list}}

\begin{document}

\title[Eigenvalue Estimates]
{Eigenvalue Estimates for Non-normal Matrices and the Zeros of
Random Orthogonal Polynomials on the Unit Circle}
\author[E.~B.~Davies and B.~Simon]{E.~B.~Davies$^{1}$ and Barry Simon$^{2}$}

\thanks{$^1$ Department of Mathematics, King's College London,
Strand, London WC2R 2LS, United Kingdom. E-mail: E.Brian.Davies@kcl.ac.uk.
Supported in part by EPSRC grant GR/R81756}
\thanks{$^2$ Mathematics 253-37, California Institute of Technology, Pasadena,
CA 91125. E-mail: bsimon@caltech.edu. Supported in part by NSF grant DMS-0140592}

\date{February 28, 2006}

\begin{abstract} We prove that for any $n\times n$ matrix, $A$, and $z$ with
$|z|\geq \|A\|$, we have that $\|(z-A)^{-1}\|\leq\cot (\frac{\pi}{4n}) \dist (z,
\spec(A))^{-1}$. We apply this result to the study of random
orthogonal polynomials on the unit circle.
\end{abstract}

\maketitle

%%%%%%%%%%%%%%%%%%%%%%%%%%%%%%%%%%%%%
\section{Introduction} \lb{s1}
%%%%%%%%%%%%%%%%%%%%%%%%%%%%%%%%%%%%%

This paper concerns a sharp bound on the approximation of eigenvalues of general
non-normal matrices that we found in a study of the zeros of orthogonal polynomials.
We begin with a brief discussion of the motivating problem, which we return to in
Section~\ref{s7}.

Given a probability measure $d\mu$ on $\bbC$ with
\begin{equation}  \lb{1.1}
\int \abs{z}^n\, d\mu(z)<\infty
\end{equation}
we define the monic orthogonal polynomials, $\Phi_n(z)$, by
\begin{gather}
\Phi_n(z) = z^n + \text{lower order}  \lb{1.2} \\
\int \ol{z^j} \, \Phi_n(z)\, d\mu(z)=0 \qquad j=0,1,\dots, n-1  \lb{1.3}
\end{gather}
If
\begin{equation}  \lb{1.4}
\begin{split}
P_n  =& \text{ orthogonal projection in $L^2(\bbC,d\mu)$} \\
&\qquad \text{onto polynomials of degree $n-1$ or less}
\end{split}
\end{equation}
then
\begin{equation}  \lb{1.5}
\Phi_n =(1-P_n)z^n
\end{equation}

A key role is played by the operator
\begin{equation}  \lb{1.6}
A_n = P_n M_z P_n \restriction \ran (P_n)
\end{equation}
where $M_z$ is the operator of multiplication by $z$ and $A_n$ is an operator
on the $n$-dimensional space $\ran(P_n)$.

If $z_0$ is a zero of $\Phi_n(z)$ of order $k$, then $f_{z_0}\equiv (z-z_0)^{-k}
\Phi_n(z)$ is in $\ran (P_n)$ and
\begin{equation}  \lb{1.7}
(A_n-z_0)^k f_{z_0}=0 \qquad (A-z_0)^{k-1} f_{z_0} \neq 0
\end{equation}
which implies
\begin{equation}  \lb{1.8}
\Phi_n (z)=\det (z-A_n)
\end{equation}
Also, $\Phi_n(z)$ is the minimal polynomial for $A_n$.

In the study of orthogonal polynomials on the real line (OPRL), a key role is
played by the fact that for any $y\in\ran(P_n)$ with $\|y\|_{L^2}=1$,
\begin{equation}  \lb{1.9}
\dist (z_0, \{\text{zeros of }\Phi_n\})\leq \|(A_n-z_0)y\|
\qquad \text{(OPRL case)}
\end{equation}
This holds because, in the OPRL case, $A_n$ is self-adjoint. Indeed, for any
normal operator, $B$, (throughout $\|\cdot\|$ is a Hilbert space norm; for
$n\times n$ matrices, the usual matrix norm induced by the Euclidean inner product)
\begin{equation}  \lb{1.10x}
\dist (z_0, \spec(B))=\|(B-z_0)^{-1}\|^{-1}
\end{equation}
and, of course, for any invertible operator $C$,
\begin{equation}  \lb{1.11x}
\inf \{\|Cy\|\mid \|y\|=1\} = \|C^{-1}\|^{-1}
\end{equation}

We were motivated by seeking a replacement of \eqref{1.9} in a case where $A_n$
is non-normal. Indeed, we had a specific situation of orthogonal polynomials
on the unit circle (OPUC; see \cite{OPUC1,OPUC2}) where one has a sequence
$z_n\in\partial\bbD = \{z\mid \abs{z}=1\}$ and corresponding unit trial vectors,
$y_n$, so that
\begin{equation}  \lb{1.9x}
\|(A_n-z_n)y_n\| \leq C_1 e^{-C_2 n}
\end{equation}
for all $n$ with $C_2 >0$. We would like to conclude that $\Phi_n(z)$
has zeros near $z_n$.

It is certainly not sufficient that $\|(A_n-z_n)y_n\|\to 0$. For the case $d\mu(z)=
d\theta/2\pi$ has $\Phi_n(z)= \dist(1,\spec(A_n))=1$, but if $y_n =(1+z+\cdots +
z^{n-1})/\sqrt{n}$, then $\|(A_n-1)y_n\|=\|P_n (z-1) y_n\|=n^{-1/2} \|P_n (z^n-1)\|
= n^{-1/2} \|1\| =n^{-1/2}$. As we will see later, by a clever choice of $y_n$,
one can even get trial vectors with $\|(A_n-1)y_n\|=O(n^{-1})$.

Of course, by \eqref{1.11x}, we are really seeking some kind of
bound relating $\|(A_n-z_n)^{-1}\|$ to $\dist (z_n, \spec(A_n))$. At
first sight, the prognosis for this does not seem hopeful. The
$n\times n$ matrix,
\begin{equation}  \lb{1.10}
N_n = \begin{pmatrix}
0 & 1 & {} & 0 \\
{} & \ddots & \ddots & {} \\
{} & {} & \ddots & 1 \\
0 & {} & {} & 0
\end{pmatrix}
\end{equation}
has
\begin{equation}  \lb{1.11}
\|(z-N_n)^{-1}\|\geq \abs{z}^{-n}
\end{equation}
since $(z-N_n)^{-1} =\sum_{j=0}^{n-1} z^{-j-1} (N_n)^j$ has $z^{-n}$ in the
$1,n$ position. Thus, as is well known, $\|(A_n-z)^{-1}\|$ for general
$n\times n$ matrices $A_n$ and general $z$ cannot be bounded by better
than $\dist (z,\spec(A_n))^{-n}$. Indeed, the existence of such bounds by
Henrici \cite{Hen} is part of an extensive literature on general variational
bounds on eigenvalues. Translated to a variational bound, this would give
$\dist (z_n,\{\text{zeros of }\Phi_n\})\leq C\|(A_n-z_n)y\|^{1/n}$, which
would not give anything useful from \eqref{1.9x}.

We note that as $n\to\infty$, there can be difficulties even if $z_0$ stays
away from $\spec(A_n)$. For, by \eqref{1.11},
\begin{equation}  \lb{1.12}
\|(1-2N_n)^{-1}\|\geq 2^{n-1}
\end{equation}
diverges as $n\to\infty$ even though $\|2N_n\|$ is bounded in $n$.

Despite these initial negative indications, we have found a linear variational
principle that lets us get information from \eqref{1.9x}. The key realization
is that $z_n$ and $\|A_n\|$ are not general. Indeed,
\begin{equation}  \lb{1.13}
\abs{z_n}=\|A_n\| =1
\end{equation}

It is not a new result that a linear bound holds in the generality we discuss.
In \cite{Nikppt}, Nikolski presents a general method for estimating norms of
inverses in terms of minimal polynomials (see the proof of Lemma~3.2 of
\cite{Nikppt}) that is related to our argument in Subsection~\ref{s6}A. His
ideas yield a linear bound but not with the optimal constant we find.

Our main theorem is

\begin{t1} Let $\calM_n$ be the set of pairs $(A,z)$ where $A$ is an
$n\times n$ matrix, $z\in\bbC$ with
\begin{equation}  \lb{1.14}
\abs{z}\geq \|A\|
\end{equation}
and
\begin{equation}  \lb{1.15}
z\notin\spec(A)
\end{equation}
Then
\begin{equation}  \lb{1.16}
c(n) \equiv \sup_{\calM_n} \, \dist (z,\spec(A)) \|(A-z)^{-1}\|  =
\cot \biggl( \f{\pi}{4n}\biggr)
\end{equation}
\end{t1}

Of course, the remarkable fact, given \eqref{1.11}, is that $c(n)<\infty$
when we only use the first power of $\dist (z,\spec(A))$. It implies that so
long as \eqref{1.14} holds,
\begin{equation}  \lb{1.16x}
\dist (z,\spec(A))\leq c(n) \|(A-z)y\|
\end{equation}
for any unit vector $y$. For this to be useful in the context of \eqref{1.9x},
we need only mild growth conditions on $c(n)$; see \eqref{1.17} below.

As an amusing aside, we note that
\begin{align*}
c(1) &=1= 0 +\sqrt1 \\
c(2) &=1 +\sqrt2 \\
c(3) &= 2+\sqrt3
\end{align*}
but the obvious extrapolation from this fails. Instead, because of properties
of $\cot(x)$,
\begin{align}
& c(n) \leq \f{4}{\pi}\, n \lb{1.17} \\
&\f{c(n)}{n}\text{ is monotone increasing to } \f{4}{\pi} \notag
\end{align}
so, in fact, for $n\geq 3$,
\[
\f{2+\sqrt3}{3}  \leq \f{c(n)}{n} \leq \f{4}{\pi}
\]
a spread of $2.3\%$.

We note that, by replacing $A$ by $A/z$ and $z$ by $1$, it suffices to prove
\begin{equation}  \lb{1.18}
\sup_{\|A\| < 1}\, \dist (1,\spec(A)) \|(1-A)^{-1}\| =
\cot \biggl(\f{\pi}{4n}\biggr)
\end{equation}
and it is this that we will establish by proving three statements. We will
use the special $n\times n$ matrix
\begin{equation}  \lb{1.19}
M_n = \begin{pmatrix}
1 & 2 & \dots & 2 \\
0 & 1 & \dots & 2 \\
\vdots & \vdots & \ddots  & \vdots \\
0 & 0 & \dots & 1
\end{pmatrix}
\end{equation}
given by
\[
(M_n)_{k\ell} = \begin{cases}
2 & \text{if } k<\ell \\
1 & \text{if } k=\ell \\
0 & \text{if } k>\ell
\end{cases}
\]
Our three sub-results are

\begin{t2} $\|M_n\|=\cot (\pi/4n)$
\end{t2}

\begin{t3} For each $0<a<1$, there exist $n\times n$ matrices $A_n(a)$ with
\begin{equation}  \lb{1.20}
\|A_n(a)\|\leq 1 \qquad \spec(A_n)=\{a\}
\end{equation}
and
\begin{equation}  \lb{1.21}
\lim_{a\uparrow 1}\, (1-a) (1-A_n (a))^{-1} = M_n
\end{equation}
\end{t3}

\begin{t4} Let $A$ be an upper triangular matrix with $\|A\| \leq 1$ and
$1\notin\spec(A)$. Then
\begin{equation}  \lb{1.22}
\dist (1, \spec(A)) \abs{(1-A)_{k\ell}^{-1}} \leq
\begin{cases}
2 & \text{if }\, k<\ell \\
1 & \text{if }\, k=\ell \\
0 & \text{if }\, k>\ell
\end{cases}
\end{equation}
\end{t4}

\begin{proof}[Proof that Theorems~2--4 $\Rightarrow$ Theorem~1] Any matrix
has an orthonormal basis in which it is upper triangular: One constructs
such a Schur basis by applying Gram-Schmidt to any algebraic basis in
which $A$ has Jordan normal form. In such a basis, \eqref{1.22} says that
\[
\dist(1, \spec(A))\|(1-A)^{-1}y\| \leq \|M_ny\| \leq \|M_n\|\, \|y\|
\]
so Theorem~2 implies LHS of \eqref{1.18} $\leq \cot(\pi/4n)$.

On the other hand, using $A_n(a)$ in $\dist (1,\spec(A))\|(1-A)^{-1}\|$
implies LHS of \eqref{1.18} $\geq \cot (\pi/4n)$. We thus have \eqref{1.18}
and, as noted, this implies \eqref{1.16}.
\end{proof}

To place Theorem~1 in context, we note that if $\abs{z}>\|A\|$,
\begin{equation}  \lb{1.23}
\|(z-A)^{-1}\| \leq \sum_{j=0}^\infty\, \abs{z}^{-j-1} \|A\|^j =
(\abs{z}-\|A\|)^{-1}
\end{equation}
So \eqref{1.16} provides a borderline between the dimension-independent
bound \eqref{1.23} for $\abs{z}>\|A\|$ and the exponential growth that
may happen if $\abs{z}<\|A\|$, essentially the phenomenon of pseudospectra
which is well documented in \cite{TrEm}; see also \cite{PG}.

The structure of this paper is as follows. In Section~\ref{s2}, we will prove
Theorem~4, the most significant result in this paper since it implies $c(n)<\infty$
and, indeed, with no effort that $c(n)\leq 2n$. Our initial proofs of $c(n)<\infty$
were more involved --- the fact that our final proof is quite simple should not
obscure the fact that $c(n)<\infty$ is a result we find both surprising and deep.

In Section~\ref{s3}, we use upper triangular Toeplitz matrices to construct
$A_n(a)$ and prove Theorem~3. Sections~\ref{s4} and \ref{s5} prove Theorem~2;
indeed, we also find that if
\begin{equation} \lb{1.27a}
(Q_n(a))_{k\ell} = \begin{cases}
1 & \text{if } k<\ell \\
a & \text{if } k=\ell \\
0 & \text{if } k > \ell
\end{cases}
\end{equation}
then
\begin{equation}  \lb{1.24}
\|Q_n(1)\| = \f{1}{2 \sin (\f{\pi}{4n+2})}
\end{equation}
which means we can compute $\|Q_n(a)\|$ for $a=0,\f12,1$. While the calculation
of $\|M_n\|$ and $\|Q_n(1)\|$ is based on explicit formulae for all the
eigenvalues and eigenvectors of certain associated operators, we could just
pull them out of a hat. Instead, in Section~\ref{s4}, we discuss the motivation
that led to our guess of eigenvectors, and in Section~\ref{s5} explicitly prove
Theorem~2.

Section~\ref{s6} contains a number of remarks and extensions concerning Theorem~1,
most importantly to numerical range concerns. Section~\ref{s7} contains the
application to random OPUC.

\medskip
\noindent{\bf Acknowledgments.} This work was done while B.~Simon was a
visitor at King's College London. He would like to thank A.~N.~Pressley and
E.~B.~Davies for the hospitality of King's College, and the London Mathematical
Society for partial support. The calculations of M.~Stoiciu \cite{StoiJAT,StoiDiss}
were an inspiration for our pursuing the estimate we found. We appreciate useful
correspondence/discussions with M.~Haase, N.~Higham, R.~Nagel, N.~K.~Nikolski,
V.~Totik, and L.~N.~Trefethen.

%%%%%%%%%%%%%%%%%%%%%%%%%%%%%%%%%%%%%
\section{The Key Bound} \lb{s2}
%%%%%%%%%%%%%%%%%%%%%%%%%%%%%%%%%%%%%

Our goal in this section is to prove Theorem~4. $A$ is an upper triangular $n\times n$
matrix. Let $\lambda_1, \dots, \lambda_n$ be its diagonal elements. Since
\begin{equation}  \lb{2.1}
\det (z-A)=\prod_{j=1}^n (z -\lambda_j)
\end{equation}
the $\lambda_j$'s are the eigenvalues of $A$ counting algebraic multiplicity.
In particular,
\begin{equation}  \lb{2.2}
\sup_j \, \abs{1-\lambda_j}^{-1} = \dist (1,\spec(A))^{-1}
\end{equation}

Define
\begin{equation}  \lb{2.3}
C=(1-A)^{-1} + (1-A^*)^{-1} -1
\end{equation}

\begin{proposition}\lb{P2.1} Suppose $\|A\|\leq 1$. Then
\begin{SL}
\item[{\rm{(a)}}]
\begin{align}
C_{jj} &= \abs{1-\lambda_j}^{-2} (1-\abs{\lambda_j}^2) \notag \\
& \leq 2\abs{1-\lambda_j}^{-1}  \lb{2.4}
\end{align}
\item[{\rm{(b)}}]
\[
C\geq 0
\]
\item[{\rm{(c)}}]
\begin{equation}  \lb{2.5}
\abs{C_{jk}} \leq \abs{C_{jj}}^{1/2} \abs{C_{kk}}^{1/2}
\end{equation}
\item[{\rm{(d)}}] If $j<k$, then $(1-A)_{jk}^{-1} = C_{jk}$.
\end{SL}
\end{proposition}

\begin{proof} (a) Since $A$ is upper triangular,
\begin{equation}  \lb{2.6}
[(1-A)^{-1}]_{jj} = (1-\lambda_j)^{-1}
\end{equation}
so \eqref{2.4} comes from
\begin{equation}  \lb{2.7}
(1-\lambda_j)^{-1} + (1-\bar\lambda_j)^{-1} -1 =
\abs{1-\lambda_j}^{-2} (1-\abs{\lambda_j}^2)
\end{equation}
and the fact that for $\abs{\lambda}\leq 1$,
\begin{align*}
\abs{1-\lambda}^{-1} (1-\abs{\lambda}^2) & = (1+\abs{\lambda})
(1-\abs{\lambda}) (\abs{1-\lambda}^{-1}) \\
&\leq 2
\end{align*}
since $1-\abs{\lambda} \leq \abs{1-\lambda}$.

\smallskip
(b) The operator analog of \eqref{2.7} is the direct computation
\begin{equation}  \lb{2.8}
C=[(1-A)^{-1}]^* (1-A^*A)(1-A)^{-1}\geq 0
\end{equation}
since $\|A\|\leq 1$ implies $A^*A\leq 1$.

\smallskip
(c) This is true for any positive definite matrix.

\smallskip
(d) $(1-A^*)^{-1}$ is lower triangular and $1$ is diagonal.
\end{proof}

\begin{proof}[Proof of Theorem~4] $(1-A)^{-1}$ is upper triangular so
$[(1-A)^{-1}]_{k\ell}=0$ if $k>\ell$. By \eqref{2.6} and \eqref{2.2},
\begin{equation}  \lb{2.9}
\abs{[(1-A)^{-1}]_{kk}} = \abs{1-\lambda_k}^{-1} \leq \dist(1,\spec(A))^{-1}
\end{equation}

By (a), (c), (d) of the proposition, if $k<\ell$,
\begin{align*}
\abs{[(1-A)^{-1}]_{k\ell}}
&\leq [\abs{1-\lambda_k}^{-2} \abs{1-\lambda_\ell}^{-2} (1-\abs{\lambda_k}^2)
(1- \abs{\lambda_\ell}^2)]^{1/2} \\
&\leq 2 [\abs{1-\lambda_k}^{-1} \abs{1-\lambda_\ell}^{-1}]^{1/2} \\
&\leq 2 [\dist (1,\spec(A))]^{-1}
\end{align*}
by \eqref{2.2}.
\end{proof}

%%%%%%%%%%%%%%%%%%%%%%%%%%%%%%%%%%%%%%%%%%%%%%%%%%%%
\section{Upper Triangular Toeplitz Matrices} \lb{s3}
%%%%%%%%%%%%%%%%%%%%%%%%%%%%%%%%%%%%%%%%%%%%%%%%%%%%

A Toeplitz matrix \cite{BSBigBk} is one that is constant along diagonals, that is,
$A_{jk}$ is a function of $j-k$. An $n\times n$ upper triangular Toeplitz matrix
(UTTM) is thus of the form
\begin{equation}  \lb{3.1}
\begin{pmatrix}
a_0 & a_1 & a_2 & \dots & a_{n-1} \\
0 & a_0 & a_1 & \dots & a_{n-2} \\
\vdots & \vdots & \vdots & \ddots & \vdots \\
0 & 0 & 0 & \cdots & a_0
\end{pmatrix}
\end{equation}
These concern us because $M_n$ is of this form and because the operators, $A_n(a)$,
of Theorem~3 will be of this form. In this section, after recalling the basics of
UTTM, we will prove Theorem~3. Then we will state some results, essentially due to
Schur \cite{Schur}, on the norms of UTTM that we will need in Section~\ref{s5} in
one calculation of the norm of $M_n$.

Given any function, $f$, which is analytic near zero, we write $T_n(f)$ for the
matrix in \eqref{3.1} if
\begin{equation}  \lb{3.2}
f(z)=a_0 + a_1 z + \cdots + a_{n-1} z^{n-1} + O(z^n)
\end{equation}
$f$ is called a symbol for $T_n(f)$.

We note that
\begin{equation}  \lb{3.3}
T_n (fg) = T_n(f) T_n(g)
\end{equation}
This can be seen by multiplying matrices and Taylor series or by manipulating
projections on $\ell^2$ (see, e.g., Corollary~6.2.3 of \cite{OPUC1}).

In addition, if $f$ is analytic in $\{z\mid \abs{z}<1\}$, then
\begin{equation}  \lb{3.4}
\|T_n(f)\| \leq \sup_{\abs{z}<1}\, \abs{f(z)}
\end{equation}
To see this well-known fact, associate an analytic function
\begin{equation}  \lb{3.5}
v(z)=v_0 + v_1 z + \cdots
\end{equation}
to the vector $\varphi_n(v)\in\bbC^n$ by
\begin{equation}  \lb{3.6}
\varphi_n(v) = (v_{n-1}, v_{n-2}, \dots, v_0)^T
\end{equation}
and note that with $\|\cdot\|_2$, the $H^2$ norm,
\begin{gather}
\|\varphi_n(v)\| = \inf \{\|v\|_2 \mid \varphi_n =\varphi_n(v)\}  \lb{3.7} \\
T_n(f) \varphi_n(v) = \varphi_n (fv) \lb{3.8}
\end{gather}
and
\begin{equation}  \lb{3.9}
\|fv\|_2 \leq \|f\|_\infty \|v\|_2
\end{equation}

If $N_n$ is given by \eqref{1.10}, then $T_n (f)=f(N_n)$, so an alternate proof of
\eqref{3.4} may be based on von Neumann's theorem; see Subsection~\ref{s6}E.

\begin{proof}[Proof of Theorem~3] For $a$ with $0<a<1$, define
\begin{equation}  \lb{3.10x}
f_a(z) = \f{z+a}{1+az}
\end{equation}
and define
\begin{equation}  \lb{3.11x}
A_n(a)=T_n (f_a)
\end{equation}

Then $f_a(e^{i\theta}) = e^{i\theta}\, \ol{(1+ae^{i\theta})}/(1+ae^{i\theta})$ has
$\abs{f_a(e^{i\theta})}=1$, so $\sup_{\abs{z}<1} \abs{f_a (z)}=1$ and thus, by
\eqref{3.4},
\begin{equation}  \lb{3.12x}
\|A_n(a)\|\leq 1
\end{equation}

By \eqref{3.1},
\begin{equation}  \lb{3.13x}
\spec(A_n(a)) = \{f_a(0)\}=\{a\}
\end{equation}

By \eqref{3.5},
\begin{equation}  \lb{3.14x}
(1-A_n(a))^{-1} = T_n ((1-f_a(z))^{-1})
\end{equation}
Now
\begin{equation}  \lb{3.15x}
(1-a)(1-f_a(z))^{-1} = \f{z+a}{1-z}
\end{equation}
so
\begin{equation}  \lb{3.10}
\lim_{a\uparrow 1}\, (1-a)(1-f_a(z))^{-1} = \f{1+z}{1-z}
\end{equation}

Thus,
\begin{equation}  \lb{3.11}
\lim_{a\uparrow 1}\, (1-a)(1-A_n(a))^{-1} = T_n \biggl( \f{1+z}{1-z}\biggr) = M_n
\end{equation}
since $(1+z)/(1-z)=1+2z+2z^2+\cdots$.
\end{proof}

We now want to refine \eqref{3.4} to get equality for a suitable $f$. A key role
is played by

\begin{lemma}\lb{L3.1} Let $\alpha\in\bbD$ and $A$ an operator with $\ol{\alpha}\,^{-1}
\notin\spec(A)$. Define
\begin{equation}  \lb{3.12}
B=(A-\alpha)(1-\ol{\alpha} A)^{-1}
\end{equation}
Then
\begin{alignat}{2}
&(1) \qquad && \|B\|\leq 1 \Leftrightarrow \|A\|\leq 1 \lb{3.13} \\
&(2) \qquad && \|B\|=1 \Leftrightarrow \|A\|=1 \lb{3.14}
\end{alignat}
\end{lemma}

\begin{proof} By a direct calculation,
\begin{equation}  \lb{3.15}
1-B^*B = (1-\alpha A^*)^{-1} [(1-\abs{\alpha}^2)(1-A^*A)](1-\ol{\alpha} A)^{-1}
\end{equation}
\eqref{3.13} follows since $1-B^*B \geq 0\Leftrightarrow 1-A^*A\geq 0$, and
\eqref{3.14} follows since \eqref{3.15} implies
\[
\inf_{\|\varphi\|=1}\, (\varphi, (1-B^*B)\varphi) = 0 \Leftrightarrow
\inf_{\|\varphi\|=1}\, (\varphi, (1-A^*A)\varphi) =0
\qedhere
\]
\end{proof}

\begin{remark} This lemma is further discussed in Subsection~\ref{s6}E.
\end{remark}

\begin{theorem}\lb{T3.2} If $A$ is an $n\times n$ UTTM with $\|A\|\leq 1$,
then there exists an analytic function, $f$, on $\bbD$ such that
\begin{equation}  \lb{3.16}
\sup_{\abs{z}<1}\, \abs{f(z)} \leq 1
\end{equation}
and
\begin{equation}  \lb{3.17}
A=T_n(f)
\end{equation}
\end{theorem}

\begin{proof} The proof is by induction on $n$. If $n=1$, $\|A\|\leq 1$ means
$\abs{a_0}\leq 1$ and we can take $f(z)\equiv a_0$. For general $n$, $\|A\|\leq 1$
means $\abs{a_0}\leq 1$. If $\abs{a_0}=1$, then $A=a_0\bdone$ and we can take
$f(z)\equiv a_0$. If $a_0 <1$, define $B$ by \eqref{3.12} with $\alpha =a_0$. $B$
is a UTTM with zero diagonal terms, so
\begin{equation}  \lb{3.17a}
B=\begin{pmatrix}
0 & {} & \ti B\\
{} & \ddots & {} \\
0 & {} & 0
\end{pmatrix}
\end{equation}
where $\|\ti B\|=\|B\|\leq 1$ by the lemma.

By the induction hypothesis, $\ti B=T_{n-1}(g)$ where
\begin{equation}  \lb{3.17b}
\sup_{\abs{z}<1}\, \abs{g(z)}\leq 1
\end{equation}
Then \eqref{3.17} holds with
\begin{equation}  \lb{3.18}
f=\f{a_0 + zg}{1+\ol{a}_0 zg}
\end{equation}
\eqref{3.17b} and \eqref{3.18} imply \eqref{3.16}.
\end{proof}

\begin{remarks} 1. By iterating $f\to g$, we see that one constructs $f$ via the
Schur algorithm; see Section~1.3 of \cite{OPUC1}.

2. Combining this and \eqref{3.4}, one obtains Schur's celebrated result
that $a_0 + a_1 z + \cdots + a_{n-1} z^{n-1}$ is the start of the Taylor
series of a Schur function if and only if the matrix $A$ of \eqref{3.1}
obeys $A^*A\leq 1$. This result is intimately connected to Nehari's theorem
on the norm of Hankel operators \cite{Neh,Pell}; see Partington \cite{Par}.

3. This is classical; see \cite{BSBigBk,NikBk,Pell}.
\end{remarks}

To state the last result of this section, we need a definition:

\begin{definition} A {\it Blaschke factor\/} is a function on $\bbD$
of the form
\begin{equation} \lb{3.19}
f(z,w) =\f{z-w}{1-\ol{w} z}
\end{equation}
where $w\in\bbD$. A (finite) {\it Blaschke product\/} is a function of
the form
\begin{equation} \lb{3.20}
f(z) = \omega\prod_{j=1}^k f(z,w_k)
\end{equation}
where $\omega\in\partial\bbD$. $k$ is called the {\it order\/} of $f$.
We allow $k=0$, in which case $f(z)$ is a constant value in $\partial\bbD$.
\end{definition}

\begin{theorem} \lb{T3.3} An $n\times n$ UTTM, $A$, has $\|A\|=c$ if and only
if $A=T_n(f)$  for an $f$ so that $c^{-1} f$ is a Blaschke product of order
$k\leq n-1$.
\end{theorem}

\begin{proof} (See as alternates: \cite{NikBk,Pell}.)
Without loss, we can take $c=1$. The proof is by induction on
$n$. If $n=1$, $k$ must be $0$, and the theorem says $\abs{a_0}=1$ if and only
if $f(0)=\omega\in\partial\bbD$, which is true.

It is not hard to see that if $f$ and $f_1$ are related by
\[
f_1(z)=z^{-1}\, \f{f(z)-f(0)}{1-\ol{f(0)}\, f(z)}
\]
then $f$ is a Blaschke product of order $k\geq 1$ if and only if $f_1$ is a
Blaschke product of order $k-1$.

Given $A$ a UTTM with $\|A\|\leq 1$, $\abs{a_0}=1$ if and only if $A=T_n(a_0)$,
that is, $A$ is given by a Blaschke product of order $0$. If $\abs{a_0}< 1$,
we define $B$ by \eqref{3.12}. $\|B\|=1$ if and only if $\|A\|=1$. $\ti B$
given by \eqref{3.17b} is related to $A$ by $A=T_n(f)$ if and only if
$\ti B=T_{n-1}(f_1)$. Thus, by induction, $\|A\|= 1$ if and only if $f$
is a Blaschke product of order $k\leq n-1$.
\end{proof}

%%%%%%%%%%%%%%%%%%%%%%%%%%%%%%%%%%%%%%%%%%%%%%%%%%%%%%%%%%%%%%%
\section{Inverse of Differential/Difference Operators} \lb{s4}
%%%%%%%%%%%%%%%%%%%%%%%%%%%%%%%%%%%%%%%%%%%%%%%%%%%%%%%%%%%%%%%

In this section and the next, we will find explicit formulae for the norms
of $M_n$ and $Q_n\equiv Q_n(1)$ given by \eqref{1.27a}. Indeed, we will find all
the eigenvalues and eigenvectors for $\abs{M_n}$ and $\abs{Q_n}$ where $\abs{A} =
\sqrt{A^*A}$. A key to our finding this was understanding a kind of continuum limit
of $M_n$: Let $K$ be the Volterra-type operator on $\calH=L^2 ([0,1],dx)$ with
integral kernel
\[
K(x,y)=\begin{cases} 1 & 0\leq x\leq y\leq 1 \\
0 & 0\leq y < x < 1
\end{cases}
\]
In some formal sense, $K$ is a limit of either $M_n$ or $Q_n$, but in a precise
sense, $M_n$ is a restriction of $K$:

\begin{proposition}\lb{P4.1} Let $\pi_n$ be the projection of $\calH$ onto the
space of functions constant on each interval $[\f{j}{n}, \f{j+1}{n})$, $j=
0,1,\dots, n-1$. Then
\begin{equation} \lb{4.1}
\pi_n K \pi_n
\end{equation}
is unitarily equivalent to $\f12 M_n/n$. In particular,
\begin{align}
\|M_n\| &\leq 2n \|K\| \lb{4.2}\\
\lim_{n\to\infty}\, \f{\|M_n\|}{n} &= 2\|K\| \lb{4.3}
\end{align}
\end{proposition}

\begin{proof} Let $\{f_j^{(n)}\}_{j=0}^{n-1}$ be the functions
\begin{equation} \lb{4.4}
f_j^{(n)}(x) = \begin{cases}
\sqrt n & \f{j}{n} \leq x < \f{j+1}{n} \\
0 & \text{otherwise} \end{cases}
\end{equation}
which form an orthonormal basis for $\ran (\pi_n)$. Since
\begin{equation} \lb{4.5}
n \langle f_j^{(n)}, K f_k^{(n)}\rangle = \tfrac12\, (M_n)_{jk}
\end{equation}
we have the claimed unitary equivalence. \eqref{4.2} is immediate from $\|\pi_n
K\pi_n\|\leq \|K\|$. \eqref{4.3} follows if we note $\slim_{n\to\infty} \pi_n=1$,
so $\lim \|\pi_n K\pi_n\|=\|K\|$.
\end{proof}

Notice that
\begin{equation} \lb{4.6}
(Kf)(x)=\int_x^1 f(y)\, dy
\end{equation}
so
\begin{equation} \lb{4.7}
\f{d}{dx}\, (Kf)=f \qquad Kf(1)=0
\end{equation}
and $K$ is an inverse of a derivative. That means $K^*K$ will be the inverse of
a second-order operator. Indeed,
\begin{align}
(K^*K)(x,y) &=\int_0^1 \ol{K(z,x)}\, K(z,y)\, dz \notag \\
&= \int_0^{\min(x,y)} dz \notag \\&= \min(x,y) \lb{4.8}
\end{align}
which, as is well known, is the integral kernel of the inverse of $-\f{d^2}{dx^2}$
with $u(0)=0$, $u'(1)=1$ boundary conditions.

We can therefore write down a complete orthonormal basis of eigenfunctions
for $K^*K$:
\begin{align}
\varphi_n(x) &= \sin(\tfrac12\, (2n-1)\pi x) \qquad n=1,2,\dots  \lb{4.9} \\
(K^*K)\varphi _n &= \f{4}{(2n-1)^2 \pi^2} \lb{4.10}
\end{align}
so
\begin{equation} \lb{4.11}
\|K\| = \|K^*K\|^{1/2} = \f{2}{\pi}
\end{equation}
By \eqref{4.2}, \eqref{4.3}, we have

\begin{corollary}
\begin{align}
\|M_n\| &\leq \f{4n}{\pi} \lb{4.12} \\
\lim_{n\to\infty}\, \f{\|M_n\|}{n} &= \f{4}{\pi} \lb{4.13}
\end{align}
\end{corollary}

Of course, we will see this when we have proven Theorem~2, but it is
interesting to have it now.

While $M_n$ is related to differential operators via \eqref{4.5}, we can
compute the norm of $Q_n$ by realizing it as the inverse of a difference
operator. Specifically, let $N_n$ be given by \eqref{1.10}. Then
\begin{equation} \lb{4.14}
(1-N_n)^{-1} = 1+N_n + N_n^2 +\cdots + N_n^{n-1}=Q_n
\end{equation}

\begin{theorem}\lb{T4.3} Let
\begin{equation} \lb{4.15}
D_n = (1-N_n)(1-N_n)^*
\end{equation}
Then $D_n$ has a complete set of eigenvectors:
\begin{align}
v_j^{(\ell)} &= \sin \biggl( \f{\pi(2\ell +1)j}{2n+1}\biggr) \qquad
j=1,\dots, n;\, \ell=0,\dots, n-1 \lb{4.16} \\
D_n v^{(\ell)} &= 4\sin^2 \biggl( \f{\pi(2\ell+1)}{2(2n+1)}\biggr) v^{(\ell)} \lb{4.17} \\
\|Q_n\| &= (\min\, \text{eigenvalue of } D_n)^{-1/2} \notag \\
&= \biggl[ 2\sin \biggl( \f{\pi}{4n+2}\biggr)\biggr]^{-1} \lb{4.18}
\end{align}
\end{theorem}

\begin{proof} By a direct calculation,
\begin{equation} \lb{4.19}
D_n = \left( \begin{array}{rrrcrrr}
2 & -1 & 0 & {} & {} & {} & {} \\
-1 & 2 & -1 & {} & {} & {} & {} \\
0 & -1 & 2 & {} & {} & {} & {} \\
{} & {} & {} & \ddots & {} & {} & {} \\
{} & {} & {} & {} & 2 & -1 & 0 \\
{} & {} & {} & {} & -1 & 2 & -1 \\
{} & {} & {} & {} & 0 & -1 & 1
\end{array} \right)
\end{equation}
is a discrete Laplacian with Dirichlet boundary condition at $0$ and Neumann
at $n$. Since
\[
-\sin (q(j+1)) + 2\sin (qj) - \sin (q(j-1)) = 4\sin^2
\biggl( \f{q}{2}\biggr) \sin (qj)
\]
\eqref{4.16}/\eqref{4.17} hold so long as $q$ is such that $\sin(q(n+1))=
\sin (qn)$, that is,
\[
\tfrac12\, [q(n+1) + qn] = (\ell + \tfrac12) \pi
\]
or $q=(2\ell +1)\pi/(2n+1)$.
\end{proof}

\begin{remark} For OPUC with $d\mu = d\theta/2\pi$, in the basis $1,z,\dots,
z^{n-1}$, $A_n$ is given by the matrix, $N_n$, of \eqref{1.10}, and so
$\|(1-N_n)^{-1}\|=\|Q_n\|\sim 2n/\pi$. Thus, there are unit vectors, $y_n$,
in this case with $\|(1-A_n)y_n\|\sim \pi/2n$.
\end{remark}

%%%%%%%%%%%%%%%%%%%%%%%%%%%%%%%%%%%%%
\section{The Norm of $M_n$} \lb{s5}
%%%%%%%%%%%%%%%%%%%%%%%%%%%%%%%%%%%%%

In this section, we will give two distinct but related proofs of Theorem~2. Both
depend on a generating function relation:

\begin{theorem}\lb{T5.1} For $\theta\in (0,\pi)$ and $z\in\bbD$, define
\begin{align}
S_\theta (z) &=\sum_{j=0}^\infty \sin ((2j+1)\theta) z^j \lb{5.1} \\
C_\theta (z) &=\sum_{j=0}^\infty \cos ((2j+1)\theta) z^j \lb{5.2}
\end{align}
Then
\begin{equation} \lb{5.3}
\f{1+z}{1-z}\, C_\theta (z) =\cot(\theta) S_\theta (z)
\end{equation}
\end{theorem}

\begin{proof} Let $\omega =e^{i\theta}$ so, summing the geometric series,
\begin{align}
S_\theta (z) &= (2i)^{-1} \sum_{j=0}^\infty (\omega^{2j+1} z^j -
\bar\omega^{2j+1} z^j) \notag \\
&= (2i)^{-1} \biggl[ \f{\omega}{1-z\omega^2}- \f{\bar\omega}{1-z\bar\omega^2}\biggr] \lb{5.4} \\
&= \f{\sin(\theta)(1+z)}{(1-z\omega^2)(1-z\bar\omega^2)} \lb{5.5}
\end{align}
For $C_\omega (z)$, the calculation is similar; in \eqref{5.4}, $(2i)^{-1}$ is replaced by
$(2)^{-1}$ and the minus sign becomes a plus:
\begin{equation} \lb{5.6}
C_\omega (z) = \f{\cos(\theta)(1-z)}{(1-z\omega^2)(1-z\bar\omega^2)}
\end{equation}
\eqref{5.5} and \eqref{5.6} imply \eqref{5.3}.
\end{proof}

Our first proof of Theorem~2 depends on looking at the Hankel matrix \cite{Par,Pell}
\begin{equation} \lb{5.7}
\wti M_n = \begin{pmatrix}
2 & 2 & \dots & 2 & 1 \\
2 & 2 & \dots & 1 & 0 \\
\vdots & \vdots & \ddots & \vdots & \vdots \\
1 & 0 & \dots & 0 & 0
\end{pmatrix}
\end{equation}
If $W_n$ is the unitary permutation matrix
\begin{equation} \lb{5.8}
(Wv)_j=v_{n+1-j}
\end{equation}
then
\begin{equation} \lb{5.9}
M_n = \wti M_n W \qquad \wti M_n =M_nW
\end{equation}
and so
\begin{equation} \lb{5.10}
\|M_n\| = \|\wti M_n\|
\end{equation}

Here is our first proof of Theorem~2:

\begin{theorem}\lb{T5.2} Let
\begin{equation} \lb{5.11}
c_j^{(n;\ell)} = \cos\biggl(\biggl(2\ell + \f12\biggr)\f{\pi}{2n}\,  (2j-1)\biggr)
\qquad j=1,2,\dots, n; \, \ell=0,\dots, n-1
\end{equation}
Then
\begin{equation} \lb{5.12}
\wti M_n c^{(n;\ell)} = \cot\biggl( \biggl( 2\ell + \f12\biggr)\f{\pi}{2n}\biggr)
c^{(n;\ell)}
\end{equation}
Thus,
\begin{equation} \lb{5.13}
\|M_n\| = \|\wti M_n\| = \cot \biggl( \f{\pi}{4n}\biggr)
\end{equation}
\end{theorem}

\begin{proof}
Let
\begin{alignat}{2}
c_j^{(n;\theta)} &=\cos (\theta(2j-1)) \qquad &&j=1,2,\dots, n \lb{5.14}  \\
\intertext{and}
s_j^{(n;\theta)} &= \sin (\theta(2j-1)) \qquad &&j=1,\dots, n \lb{5.15}
\end{alignat}
Then \eqref{5.3} implies that
\begin{equation} \lb{5.16}
M_n W c^{(n;\theta)} = \cot(\theta) W s^{(n;\theta)}
\end{equation}
by looking at coefficients of $1,z,\dots, z^{n-1}$. The $W$ comes from
\eqref{3.6}/\eqref{3.8}. If
\begin{equation} \lb{5.17}
\theta = \f{\pi}{2} + 2\ell\pi \qquad \ell=0, \dots, n-1
\end{equation}
then
\begin{equation} \lb{5.18}
W s^{(n;\theta)} = c^{(n;\theta)}
\end{equation}
and \eqref{5.16} becomes \eqref{5.12}.

Since $\wti M$ is self-adjoint, \eqref{5.13} follows from \eqref{5.12} either by noting
that $\max \abs{\cot ((2\ell + \f12) \f{\pi}{2n})}=\cot (\f{\pi}{4n})$ or by noting that
$c^{(n;\theta = \pi/4n)}$ is a positive eigenvector of a positive self-adjoint matrix,
so its eigenvalue is
the norm by the Perron-Frobenius theorem.
\end{proof}

Our second proof relies on the following known result (see Milovani\'c et al.\
\cite{MMR}, page~272, and references therein; this result is called the
Enestr\"om-Kakeya theorem; see also P\'olya-Szeg\H{o} \cite{PS}, problem 22
on pp.~107 and 301, who also mention Hurwitz):

\begin{lemma} \lb{L5.3} Suppose
\begin{equation} \lb{5.18a}
0 < a_0 < a_1 < \cdots < a_n
\end{equation}
Then
\begin{equation} \lb{5.19}
P(z)=a_0 + a_1 z+\cdots + a_n z^n
\end{equation}
has all its zeros in $\bbD$.
\end{lemma}

\begin{theorem}\lb{T5.4} Let
\begin{align}
S^{(n)}(z) &=\sum_{j=0}^{n-1} \sin \biggl( (2j+1)\f{\pi}{4n}\biggr) z^j \lb{5.20} \\
C^{(n)}(z) &=\sum_{j=0}^{n-1} \cos\biggl( (2j+1) \f{\pi}{4n}\biggr) z^j \lb{5.21}
\end{align}
Then
\begin{equation} \lb{5.22}
b^{(n)}(z) = \f{S^{(n)}(z)}{C^{(n)}(z)}
\end{equation}
is a Blaschke product of order $n-1$. Moreover,
\begin{equation} \lb{5.23}
\cot\biggl( \f{\pi}{4n}\biggr) b^n (z) = 1 +2\sum_{j=1}^{n-1} z^j + O(z^n)
\end{equation}
and
\begin{equation} \lb{5.24}
\|M_n\| = \cot\biggl( \f{\pi}{4n}\biggr)
\end{equation}
\end{theorem}

\begin{proof} The coefficients of $S^{(n)}$ obey \eqref{5.18a} so, by the lemma,
$S^{(n)}$ has all its zeros in $\bbD$. Moreover, by \eqref{5.18}, $C^{(n)}(z)=z^n \,
\ol{S^{(n)}(1/\bar z)}$, which implies \eqref{5.22} is a Blaschke product.

\eqref{5.23} is just a translation of \eqref{5.3}. \eqref{5.23} implies \eqref{5.24}
by Theorem~\ref{T3.3}.
\end{proof}

%%%%%%%%%%%%%%%%%%%%%%%%%%%%%%%%%%%%%%%%%%%%%
\section{Some Remarks and Extensions} \lb{s6}
%%%%%%%%%%%%%%%%%%%%%%%%%%%%%%%%%%%%%%%%%%%%%

In this section,we make some remarks that shed light on or extend Theorem~1, our main
result.

\subsection*{A. An alternate proof}
We give a simple proof of a weakened version of Theorem~4 but which suffices for
applications like those in Section~\ref{s7}. This argument is related to ones in
Section~3 of Nikolski \cite{Nikppt}.

\begin{theorem} \lb{T6.1} If $\|A\|\leq 1$ and $1\notin\spec(A)$, then
\begin{equation} \lb{6.1}
\dist(1,\spec(A)) \|(1-A)^{-1}\|  \leq 2m
\end{equation}
where $m$ is the degree of the minimal polynomial for $A$.
\end{theorem}

\begin{proof} We prove the result for $\|A\| < 1$. The general result follows by
taking limits. We make repeated use of Lemma~\ref{L3.1} which implies that if,
for $\lambda\in\bbD$, and we define
\begin{equation} \lb{6.2}
B(\lambda) = \biggl(\f{A-\lambda}{1-\ol{\lambda} A}\biggr)
\biggl( \f{1-\ol{\lambda}}{1-\lambda}\biggr)
\end{equation}
then
\begin{equation} \lb{6.3}
\|B(\lambda)\|\leq 1
\end{equation}
By algebra,
\begin{equation} \lb{6.4}
(1-x)^{-1} \biggl[ 1-\f{x-\lambda}{1-\ol{\lambda} x}
\biggl( \f{1-\ol{\lambda}}{1-\lambda}\biggr)\biggr]
= \f{1}{1-\lambda} \, \biggl[ 1+ \ol{\lambda}
\biggl( \f{x-\lambda}{1-x \ol{\lambda}}\biggr)\biggr]
\end{equation}
so, by Lemma~\ref{L3.1} again,
\begin{equation} \lb{6.5}
\|(1-A)^{-1} (1-B(\lambda))\| \leq \abs{1-\lambda}^{-1} (1+\abs{\lambda})
\end{equation}

Now let $\prod_{j=1}^m (x-\lambda_j)$ be the minimal polynomial for $A$. Then
\[
\prod_{j=1}^m B(\lambda_j) =0
\]
so
\begin{align}
(1-A)^{-1} &= (1-A)^{-1} \biggl[ 1-\prod_{j=1}^m B_j (\lambda) \biggr] \notag \\
&=\sum_{j=1}^m\, (1-A)^{-1} [1- B_j(\lambda)] \prod_{k=j+1}^m B_k (\lambda) \lb{6.6}
\end{align}
(the empty product for $j=m$ is interpreted as the identity operator) which, by
\eqref{6.3} and \eqref{6.5}, implies
\begin{align*}
\text{LHS of \eqref{6.1}} &\leq \sum_{j=1}^m \dist (1,\spec(A)) \abs{1-\lambda_j}^{-1}
(1+\abs{\lambda_j}) \\
&\leq 2m
\end{align*}
since $1+\abs{\lambda_j}\leq 2$ and $\lambda_j\in\spec(A)$ so $\dist (1,\spec(A))
\abs{1-\lambda_j}^{-1}\leq 1$.
\end{proof}

\begin{remarks} 1. The factor $(1- \ol{\lambda})/(1-\lambda)$ is taken in \eqref{6.2} so
$f_\lambda(z) = (z-\lambda)(1- \ol{\lambda} z)^{-1} (1- \ol{\lambda})(1-\lambda)^{-1}$ has
$1-f_\lambda (1)=0$.

2. In place of the algebra \eqref{6.4}, one can compute that the $\sup_{\abs{z}<1}
\text{LHS of \eqref{6.4}}$ is $\abs{1-\lambda}^{-1} [1+ \abs{\lambda}]$ and use
von Neumann's theorem as discussed in Subsection~E below.
\end{remarks}

\subsection*{B. Minimal polynomials}

While the constant $2$ in \eqref{6.1} is worse than $4/\pi$ in \eqref{1.16}/\eqref{1.17},
\eqref{6.1} appears to be stronger in that $m$, not $n$, appears, but we can also
strengthen \eqref{1.16} in this way:

\begin{theorem}\lb{T6.2} If $\|A\|\leq 1$, $1\notin\spec(A)$, and $m$ is the degree of
the minimal polynomial for $A$, then
\begin{equation} \lb{6.7}
\dist(1,\spec(A)) \|(1-A)^{-1}\| \leq \cot \biggl( \f{\pi}{4m}\biggr)
\end{equation}
\end{theorem}

\begin{proof} Let $\|y\|=1$. Since $A^m y$ is a linear combination of
$\{A^j y\}_{j=0}^{m-1}$, the cyclic subspace, $V_y$, has $\dim (V_y)
\equiv m_y\leq m$. Since $A\restriction V_y$ is an operator of a space
of dimension $m_y$, we have
\begin{align*}
\dist(1,\spec(A))\|(1-A)^{-1}y\| &\leq c(m_y) = \cot\biggl( \f{\pi}{4m_y}\biggr) \\
&\leq \cot\biggl( \f{\pi}{4m}\biggr)
\qedhere
\end{align*}
\end{proof}

\subsection*{C. Numerical range}

For any bounded operator, $A$, on a Hilbert space, the numerical range, $\num(A)$, is defined by
\begin{equation} \lb{6.8}
\num(A)=\{\langle \varphi, A\varphi\rangle\mid \|\varphi\|=1\}
\end{equation}
It is a bounded convex set (see \cite[p.~150]{Dav80}), and when $A$ is a finite matrix,
also closed. Theorem~1 can be improved to read:

\begin{theorem}\lb{T6.3} Let $\wti\calM_n$ be the set of pairs $(A,z)$ where $A$
is an $n\times n$ matrix, $z\in\bbC$ with
\begin{equation} \lb{6.9}
z\notin\spec(A) \qquad z\notin\num(A)^\intt
\end{equation}
Then
\begin{equation} \lb{6.10}
\sup_{\wti\calM_n}\, \dist (z, \spec(A)) \|(A-z)^{-1}\| =
\cot \biggl( \f{\pi}{4n}\biggr)
\end{equation}
\end{theorem}

\begin{remarks} 1. Since $\num(A)\subset \{z\mid \abs{z}\leq \|A\|\}$,
$\calM_n\subset \wti\calM_n$, and this is a strict improvement of
\eqref{1.16}.

2. We need only prove
\[
\dist(z,\spec(A)) \|(A-z)^{-1}\|\leq \cot\biggl( \f{\pi}{4n}\biggr)
\]
since the equality then follows from $\calM_n\subset\wti\calM_n$.

3. By replacing $A$ by $e^{i\theta}(A-z)$ for suitable $\theta$ and $z$, we
need only prove
\begin{equation} \lb{6.11}
\Real (A) \geq 0, \, 0\notin\spec(A) \Rightarrow \dist (0,\spec(A))\|A^{-1}\|
\leq \cot\biggl( \f{\pi}{4n}\biggr)
\end{equation}
for by convexity of $\num(A)$, if $z\notin\num(A)^\intt$, there is a half-plane,
$P$, with $\num(A)\subset P$ and $z\in\partial P$. It is \eqref{6.11} we will
prove below.
\end{remarks}

\begin{proof}[First Proof of Theorem~\ref{T6.3}] Let
\begin{align}
C &= A^{-1} + (A^*)^{-1} \lb{6.12} \\
&= (A^*)^{-1} 2\Real (A)(A)^{-1} \geq 0 \lb{6.13}
\end{align}
Thus,
\begin{equation} \lb{6.14}
\abs{C_{jk}} \leq \abs{C_{jj}}^{1/2} \abs{C_{kk}}^{1/2}
\end{equation}
Now just follow the proof of Theorem~4 in Section~\ref{s2}.
\end{proof}

\begin{proof}[Second Proof of Theorem~\ref{T6.3}] We use Cayley transforms.
For $0<s$, define
\begin{equation} \lb{6.15}
B(s) =(1-sA)(1+sA)^{-1}
\end{equation}
Since
\[
\|(1+sA)\varphi\|^2 -\|(1-sA)\varphi\|^2 = 4s \, \Real (\varphi, A\varphi)\geq 0
\]
we have that
\begin{equation} \lb{6.16}
\|B(s)\|\leq 1
\end{equation}

Because
\begin{equation} \lb{6.17}
1-B(s) =2sA (1+sA)^{-1}
\end{equation}
we have for $s$ small that
\begin{equation} \lb{6.18}
\dist(1,\spec(B(s))) = 2s\, \dist(0,\spec(A)) + O(s^2)
\end{equation}
Thus, by Theorem~1,
\begin{equation} \lb{6.19}
2s \, \dist(0,\spec(A)) \|(1-B(s))^{-1}\| \leq \cot\biggl( \f{\pi}{4n}\biggr)
+ O(s)
\end{equation}

By \eqref{6.17},
\[
(1-B(s))^{-1} = (2s)^{-1} [A^{-1} +s]
\]
so
\begin{equation} \lb{6.20}
\|A^{-1}\| \leq \abs{s} + 2s \|(1-B(s))^{-1}\|
\end{equation}
This plus \eqref{6.18} implies \eqref{6.11} as $s\downarrow 0$.
\end{proof}

\subsection*{D. Bounded powers}

We note that there is also a result if
\begin{equation} \lb{6.21}
\sup_{m\geq 0}\, \|A^m\| =c < \infty
\end{equation}
We suspect the $3/2$ power in the following is not optimal. We note that one can
also use this method if $\|A^m\|$ is polynomially bounded in $m$.

\begin{theorem}\lb{T6.4} If \eqref{6.21} holds, then
\begin{equation} \lb{6.22}
\|(1-A)^{-1}\| \leq c(3n)^{3/2} \dist (1,\spec(A))^{-3/2}
\end{equation}
\end{theorem}

\begin{proof} By the argument of Section~\ref{s1} (using \eqref{1.11x}), this is
equivalent to
\begin{equation} \lb{6.23}
\dist(1,\spec(A)) \leq 3n (c\|(1-A)y\|)^{2/3}
\end{equation}
for all unit vectors $y$.

Define for $1<r$,
\begin{equation} \lb{6.24}
\langle f,g\rangle_r = \sum_{m=0}^\infty r^{-2m} \langle A^m f, A^m g\rangle
\end{equation}
By \eqref{6.21},
\begin{equation} \lb{6.25}
\|f\| \leq \|f\|_r \leq cr (r^2-1)^{-1/2} \|f\|
\end{equation}

By \eqref{6.24},
\begin{equation} \lb{6.26}
\|Af\|_r^2 \leq r^2 \|f\|_r^2
\end{equation}
so
\begin{equation} \lb{6.27}
\|A\|_r \leq r
\end{equation}
so if $C=r^{-1}A$, then
\begin{equation} \lb{6.28}
\|C\|_r \leq 1
\end{equation}

Clearly, for $\|y\|=1\leq \|y\|_r$,
\begin{align}
\|Cy-y\|_r &\leq \abs{r^{-1}-1}\, \|y\|_r + r^{-1} \|(A-1)y\|_r \notag \\
&\leq \abs{r^{-1}-1}\, \|y\|_r + c (r^2 -1)^{-1/2} \|(A-1)y\| \notag \\
&\leq ((r-1) + c[2(r-1)]^{-1/2} \|(A-1)y\|) \|y\|_r \lb{6.29}
\end{align}

It follows by Theorem~1 and the fact that $\spec(A)$ is independent of
$\|\cdot\|_r$ that
\begin{equation} \lb{6.30}
\dist(1, r^{-1}\spec(A)) \leq \f{4n}{\pi}\, \{c\|(A-1)y\| (2(r-1))^{-1/2} +
(r-1)\}
\end{equation}
and thus
\begin{equation} \lb{6.31}
\dist(1,\spec(A)) \leq (r-1) + \f{4\pi}{n}\, \{c\|(A-1)y\| (2(r-1))^{-1/2} +
(r-1)\}
\end{equation}

Choosing $r=1 + \f12 (c\|(A-1)y\|)^{2/3}$ and using $\f12 + \f{6n}{\pi}\leq 3n$, we
obtain \eqref{6.23}.
\end{proof}

\subsection*{E. Von Neumann's theorem}

Lemma~\ref{L3.1} is a special case of a theorem of von Neumann. The now standard
proof of this result uses Nagy dilations \cite{FSz}; we have found a simple
alternative that relies on

\begin{lemma}\lb{L6.5} For any $A$, with $\|A\| < 1$ and $A=U\abs{A}$, and $U$
unitary, there exists an operator-valued function, $g$, analytic in a neighborhood
of $\ol{\bbD}$ so that $g(e^{i\theta})$ is unitary and $g(0)=A$.
\end{lemma}

\begin{proof} Let
\begin{equation} \lb{6.35}
g(z)=U\biggl[ \f{z+\abs{A}}{1+z\abs{A}}\biggr]
\end{equation}

The factor in $[\dots]$ is unitary if $z=e^{i\theta}$, since
\begin{align*}
(e^{i\theta}+\abs{A})^* (e^{i\theta} + \abs{A})
&= 1 + A^* A + 2\cos \theta \abs{A} \\
&= (1+e^{i\theta} \abs{A})^* (1+e^{i\theta} \abs{A})
\qedhere
\end{align*}
\end{proof}

\begin{theorem}[von Neumann \cite{vNeu}] \lb{T6.6A} Let $f\colon \bbD\to\bbD$. If $\|A\|<1$,
define $f(A)$ by
\begin{equation} \lb{6.32}
f(z) =\sum_{n=0}^\infty a_n z^n \qquad
f(A) \equiv \sum_{n=0}^\infty a_n A^n
\end{equation}
Then
\begin{equation} \lb{6.33}
\|f(A)\| \leq 1
\end{equation}
\end{theorem}

\begin{proof}[Proof of von Neumann's theorem, given the lemma] \ Suppose first
that $A$ obeys the hypotheses of the lemma. By a limiting argument, suppose $f$
is analytic in a neighborhood of $\ol{\bbD}$. Applying the maximum principle to
$f(g(z))$, we see
\begin{align}
\|f(A)\| &= \|f(g(0))\| \leq \sup_\theta \, \|f(g(e^{i\theta}))\| \notag \\
&= \sup_\theta\, \abs{f(e^{i\theta})} \leq 1 \lb{6.34}
\end{align}
where \eqref{6.34} uses the spectral theorem for the unitary $g(e^{i\theta})$.

For general $A$, if $\ti A=A\oplus 0$ on $\calH\oplus\calH$, then $\ti A=U\abs{\ti A}$
with $U$ unitary and we obtain $\|f(\ti A)\|\leq 1$. But $f(\ti A)=f(A)\oplus 0$.
\end{proof}

\begin{remarks} 1. In general, $A=V\abs{A}$ with $V$ a partial isometry. We can extend
this to a unitary $U$ so long as $\dim (\ran (V)^\perp)=\dim (\ker (V)^\perp)$.
This is automatic in the finite-dimensional case and also if $\dim(\calH)=
\infty$ for $A\oplus 0$ since then both spaces are infinite-dimensional.

2. This proof is close to one of Nelson \cite{Nel} who also uses the maximum
principle and polar decomposition, but uses a different method for interpolating
the self-adjoint part (see also Nikolski \cite{NikBk}).
\end{remarks}

%%%%%%%%%%%%%%%%%%%%%%%%%%%%%%%%%%%%%
\section{Zeros of Random OPUC} \lb{s7}
%%%%%%%%%%%%%%%%%%%%%%%%%%%%%%%%%%%%%

In this section, we apply Theorem~1 to obtain results on certain OPUC. We begin
by recalling the recursion relations for OPUC \cite{OPUC1,OPUC2,1Ft}. For each
non-trivial probability measure, $d\mu$, on $\partial\bbD$, there is a sequence
of complex numbers, $\{\alpha_n(d\mu)\}_{n=0}^\infty$, called Verblunsky
coefficients so that
\begin{equation} \lb{7.1}
\Phi_{n+1}(z)=z\Phi_n(z) -\bar\alpha_n \Phi_n^*(z)
\end{equation}
where
\begin{equation} \lb{7.2}
\Phi_n^*(z) = z^n\, \ol{\Phi_n (1/\bar z)}
\end{equation}

The $\alpha_n$ obey $\abs{\alpha_n}<1$ and Verblunsky's theorem \cite{OPUC1,1Ft}
says that $\mu\mapsto\{\alpha_n(d\mu)\}_{n=0}^\infty$ is a bicontinuous bijection
from the non-trivial measures on $\partial\bbD$ with the topology of vague
convergence to $\bbD^\infty$ with the product topology.

For each $\rho$ in $(0,1)$, we define the $\rho$-model to be the set of
random Verblunsky coefficients where $\alpha_n$ are independent,
identically distributed random variables, each uniformly distributed in
$\{z\mid\abs{z}\leq\rho\}$. A point in the model space of $\alpha$'s will
be denoted $\omega$; $\Phi_n (z;\omega)$ will be the corresponding OPUC and
$\{z_j^{(n)}(\omega)\}_{j=1}^n$ the zeros of $\Phi_n$ counting multiplicity.
Our results here depend heavily on earlier results of Stoiciu \cite{StoiJAT,StoiDiss},
who studied a closely related problem (see below). In turn, Stoiciu relied, in part,
on earlier work on eigenvalues of random Schr\"odinger operators \cite{Mol,Min}.

We will prove the following three theorems:
\begin{theorem}\lb{T7.1} Let $0<\rho <1$. Let $k\in \{1,2,\dots\}$. Then for
a.e.\ $\omega$ in the $\rho$-model,
\begin{equation} \lb{7.3}
\limsup_{n\to\infty}\, \f{\# \{ j\mid \abs{z_j^{(n)}(\omega)} < 1-n^{-k}\}}{[\log(n)]^2}
<\infty
\end{equation}
\end{theorem}

Thus, the overwhelming bulk of zeros are polynomially close to $\partial\bbD$.
If we look at a small slice of argument, we can say more:

\begin{theorem}\lb{T7.2} Let $0<\rho < 1$. Let $\theta_0\in [0,2\pi)$ and $a<b$
real. Let $\eta <1$. Then with probability $1$, for large $n$, there are no zeros
in $\{z\mid \arg z\in (\theta_0 + \f{2\pi a}{n}, \theta_0 + \f{2\pi b}{n});\,
\abs{z}<1 -\exp (-n^\eta)\}$.
\end{theorem}

Finally and most importantly, we can describe the statistical distribution
of the arguments:

\begin{theorem}\lb{T7.3} Let $0<\rho <1$. Let $\theta_0\in [0,2\pi)$. Let
$a_1 < b_1 \leq a_2 < b_2 \leq \cdots \leq a_\ell <b_\ell$ and let $k_1,
\dots, k_\ell$ be in $\{0,1,2,\dots\}$. Then as $n\to\infty$,
\begin{equation} \lb{7.4}
\prob \biggl( \# \biggl( j\biggm| \arg z_j^{(n)}(\omega)\in
\biggl( \theta_0 + \f{2\pi a_m}{n}\, , \theta_0 + \f{2\pi b_n}{n}\biggr)\biggr)
= k_m \text{ for } m=1, \dots, \ell
\end{equation}
converges to
\begin{equation} \lb{7.5}
\prod_{m=1}^\ell \f{(b_m-a_m)^{k_m}}{k_m!}\, e^{-(b_m -a_m)}
\end{equation}
\end{theorem}

This says the zeros are asymptotically Poisson distributed. As we stated,
our proofs rely on ideas of Stoiciu, essentially using Theorem~1 to complete
his program. To state the results of his that we use, we need a definition.

For $\beta\in\partial\bbD$, the paraorthogonal polynomials (POPUC) are
defined by
\begin{equation} \lb{7.6}
\Phi_n^{(\beta)}(z)=\Phi_{n-1}(z) - \bar\beta \Phi_{n-1}^*(z)
\end{equation}
These have zeros on $\partial\bbD$. Indeed, they are eigenvalues of a rank
one unitary perturbation of the operator $A_n$ of \eqref{1.6}. We extend
the $\rho$-model to include an additional set of independent parameters
$\{\beta_j\}_{j=0}^\infty$ in $\partial\bbD$, each uniformly distributed on
$\partial\bbD$. $\ti z_j^{(n)}(\omega)$ denotes the zeros of $\Phi_n^{(\beta_n)}
(z;\omega)$. Stoiciu \cite{StoiJAT,StoiDiss} completely analyzed these POPUC
zeros. We will need three of his results:

\begin{theorem}[$=$ Theorem~6.1.3 of \cite{StoiDiss} $=$ Theorem~6.3 of \cite{StoiJAT}]
\lb{T7.4} Let $I$ be an interval in $\partial\bbD$. Then
\begin{equation} \lb{7.7}
\prob(2 \text{ or more $\ti z_j^{(n)}(\omega)$ lie in $I$})\leq
\f12\biggl( \f{n\abs{I}}{2\pi}\biggr)^2
\end{equation}
where $\abs{I}$ is the $d\theta$ measure of $I$.
\end{theorem}

For the next theorem, we need the fact that there is an explicit realization of
$A_n$ and the associated rank one perturbations as $n\times n$ complex CMV matrices
(see \cite{CMV,OPUC1,OPUC2,1Ft}), $\calC_n$, whose eigenvalues are the $z_j^n$, and
$\ti\calC_n^{(\beta_n)}$ whose eigenvalues are the $\ti z_j^n$, so that
\begin{equation} \lb{7.8}
\|(\calC_n - \calC_n^{(\beta_n)})\varphi\| \leq \abs{\varphi_{n-1}} + \abs{\varphi_n}
\end{equation}

The next theorem uses the components so \eqref{7.8} holds.

\begin{theorem}[$=$ Theorem~1.1.2 of \cite{StoiDiss} $=$ Theorem~2.2 of \cite{StoiJAT}]
\lb{T7.5} There exists a constant $D_2$ {\rm{(}}depending only on $\rho${\rm{)}} so
that for every eigenvector $\varphi^{(j,\omega; n)}$ of $\ti\calC_n^{(\beta_n)}$,
we have for
\begin{equation} \lb{7.8a}
\abs{m-m(\varphi^{(j,\omega;n)})}\geq D_2 (\log n)
\end{equation}
that
\begin{equation} \lb{7.9}
\abs{\varphi_m^{(j,\omega;n)}} \leq C_\omega e^{-4\abs{m-m(\varphi^{(j,\omega;n)})}/D_2}
\end{equation}
where $C_\omega$ is an a.e.\ finite constant and
\begin{equation} \lb{7.9a}
m(\varphi) = \text{first $k$ so } \abs{\varphi_k} =
\max_m\, \abs{\varphi_m}
\end{equation}
\end{theorem}

We will also need the results that Stoiciu proves along the way that for each $C_0$,
\begin{equation} \lb{7.9b}
\{\omega\mid C_\omega < C_0\} \equiv \Omega_{C_0}
\end{equation}
is invariant under rotation of the measures $d\mu_\omega$, and that for each $C_0$ fixed
and all $\omega\in\Omega_{C_0}$,
\begin{equation} \lb{7.10}
\#(j\mid m(\varphi^{(j,\omega;n)}) =m_0) \leq D_3 (\log n)
\end{equation}
where $D_3$ is only $C_0$-dependent and is independent of $\omega$, $m_0$, and $n$.
\eqref{7.10} comes from the fact that, by \eqref{7.9}, for $D_3$ only depending on $C_0$,
\begin{equation} \lb{7.11}
\sum_{\abs{m-m(\varphi)}\geq  \f14 D_3(\log n)}\,\abs{\varphi_m}^2 \leq \tfrac12
\end{equation}
so, by \eqref{7.9a}, for $\varphi$'s with $m(\varphi)=m_0$,
\begin{equation} \lb{7.12}
\tfrac12\, D_3 (\log n) \abs{\varphi_{m_0}}^2 \geq \tfrac12
\end{equation}
which, given
\begin{equation} \lb{7.13}
\sum_\varphi\, \abs{\varphi_{m_0}}^2 =1
\end{equation}
implies \eqref{7.10}.

The last of Stoiciu's results we will need is

\begin{theorem}[$=$ Theorem~1.0.6 of \cite{StoiDiss} $=$ Theorem~1.1 of \cite{StoiJAT}]
\lb{T7.6} For $\theta_0\in [0,2\pi)$ and $a_1<b_1 \leq a_2 < b_2\leq \cdots\leq a_\ell <
b_\ell$ and $k_1, \dots, k_\ell$ in $\{0,1,2,\dots\}$, we have, as $n\to\infty$, that
\eqref{7.4} with $z_j^{(n)}$ replaced by $\ti z_j^{(n)}$ converges to \eqref{7.5}.
\end{theorem}

With this background out of the way, we begin the proofs of the new
Theorems~\ref{T7.1}--\ref{T7.3} with

\begin{theorem}\lb{T7.7} Fix $\rho\in (0,1)$. Then for a.e.\ $\omega$, there exists
$N_\omega$ so if $n\geq N_\omega$, then
\begin{equation} \lb{7.14}
\min_{j\neq k}\, \abs{\ti z_j^{(n)} - \ti z_k^{(n)}} \geq 2n^{-4}
\end{equation}
\end{theorem}

\begin{remark} $n^{-3-\veps}$ will work in place of $n^{-4}$.
\end{remark}

\begin{proof} For each $n$, cover $\partial\bbD$ by two sets of intervals of size
$4n^{-4}$: one set non-overlapping, except at the end, starting with $[0, 4n^{-4}]$
and the other set starting with $[2n^{-4}, 6n^{-4}]$. If \eqref{7.14} fails for some
$n$, then there are two zeros within one of these intervals. By \eqref{7.7}, the
probability of two zeros in one of these intervals is $O((nn^{-4})^2)=O(n^{-6})$.
The number of intervals at order $n$ is $O(n^4)$. Since $\sum_{n=1}^\infty n^4 n^{-6}
<\infty$, the sum of the probabilities of two zeros in an interval is summable. By the
Borel-Cantelli lemma \cite{Stroo} for a.e.\ $\omega$, only finitely many intervals
have two zeros. Hence, for large $n$, \eqref{7.14} holds.
\end{proof}

\begin{proof}[Proof of Theorem~\ref{T7.1}] Obviously, if \eqref{7.3} holds for some $k$,
it holds for all smaller $k$, so we will prove it for $k\ge 4$. We also need only prove
it on any $\Omega_{C_0}$ given by \eqref{7.9b} since $\cup\Omega_{C_0}$ has probability
$1$ by Theorem~\ref{T7.5}. Consider those $\varphi^{(j,\omega;n)}$ with
\begin{equation} \lb{7.15}
\abs{m(\varphi^{(j,\omega;n)})-n} \geq K (\log n)
\end{equation}
By \eqref{7.10}, the number of $j$ for which \eqref{7.15} fails is $O((\log n)^2)$.

By \eqref{7.9} and \eqref{7.8} and the fact that $\varphi$ is a unit eigenfunction,
then
\begin{equation} \lb{7.16}
\|(\calC_n -\ti z_j^{(n)})\varphi^{(j,\omega;n)}\|\leq 2 C_\omega n^{-4K/D_2}
\end{equation}
so picking $K$ large enough and $n$ large enough that $\f{4}{\pi} 2C_\omega n^{-1}
<1$, we have
\begin{equation} \lb{7.17}
\|(\calC_n -\ti z_j^{(n)})\varphi^{(j,\omega;n)}\| \leq\f{\pi}{4n}\, n^{-k}
\end{equation}

Thus, by Theorem~1 and $\|\calC_n\|=1 = \abs{\ti z_j^{(n)}}$, we see that for
each $j$ obeying \eqref{7.15}, there is a $z_j^{(n)}$ so
\begin{equation} \lb{7.18}
\abs{z_j^{(n)}-\ti z_j^{(n)}}\leq n^{-k}
\end{equation}
By Theorem~\ref{T7.7} and $k\geq 4$, the $z_j^{(n)}$ are distinct for $n$ large,
so we have $n-O((\log n)^2)$ zeros with $\abs{z_j^{(n)}}\geq 1-n^{-k}$. This is
\eqref{7.3}.
\end{proof}

\begin{proof}[Proof of Theorem~\ref{T7.2}] In place of \eqref{7.15}, we look for
$\varphi$'s so
\begin{equation} \lb{7.19}
\abs{m(\varphi^{(j,\omega;n)})-n}\geq \f{D_2}{2}\, n^{1-\eta}
\end{equation}
For such $j$'s, using the above arguments, there are zeros $z_j^{(n)}$ with
\begin{equation} \lb{7.20}
\abs{z_j^{(n)} -\ti z_j^{(n)}} \leq C_\omega \exp (-2n^\eta)
\end{equation}
\end{proof}

As in Stoiciu \cite{StoiJAT,StoiDiss}, the distribution of $\ti z_j^{(n)}$
for which \eqref{7.19} fails is rotation invariant. Since the number is $O(n^{1-\eta}
\log n)$ out of $O(n)$ zeros, the probability of any of these had zeros lying in
$\{z\mid\arg z\in (\theta_0 + \f{2\pi a}{n}, \theta_0 + \f{2\pi b}{n})\}$ goes to
zero as $n\to\infty$.

\begin{proof}[Proof of Theorem~\ref{T7.3}] By the last proof, the zeros of $\Phi_n$
with the given arguments lie within $O(e^{-n^\eta})$ of those of $\Phi_n^{(\beta)}$
and, by Theorem~\ref{T7.7}, these zeros are distinct. Theorem~\ref{T7.6} completes
the proof if one gets upper and lower bounds by slightly increasing/decreasing the
intervals on an $O(1/n)$ scale.
\end{proof}

We close with the remark about improving these theorems. While \eqref{7.10} is the best
one can hope for as a uniform bound, with overwhelming probability the number should be
bounded. Thus, we expect in Theorem~\ref{T7.1} that one can obtain $O((\log n)^{-1})$
in place of $O((\log n)^{-2})$. It is possible in Theorem~\ref{T7.2} that one can
improve $O(e^{-n^\eta})$ for all $\eta\in 1$ to $O(e^{-An})$ for some $A$.

\bigskip
%%%%%%%%%%%%%%%%%%%%%%%%%%%%%

\end{document}